\documentclass{article}
\usepackage{amsmath}
\newtheorem{thm}{Theorem}
\usepackage{graphicx}
\title{Applying Adaptive Gradient Descent to solve matrix factorization}
\author{Dan Qiao}
\begin{document}
\maketitle
\begin{abstract}
	This note is based on [1], which studies the minimization of a convex function f(X) over the set of $n\times n$ positive semi-definite matrices. Similar to [1], we optimize the function $g(U)=f(UU^{T})$, with $U\in R_{n\times r}$ and $r\leq n$. This is the method called Factored gradient descent(FGD). We apply the method of adaptive gradient descent which uses different step length at different epoch. Adaptive gradient descent performs much better than FGD in the tests and keeps the guarantee of convergence speed at the same time.
\end{abstract}
	\section{Introduction}
	Consider the following standard convex semi-definite optimization problem: 
	$$minimize_{X\in R^{n\times n}} f(X) \ subject\ to\ X \geq 0,$$ where f: $R^{n\times n}\rightarrow R$ is a convex and differentiable function, and $X\geq 0$ denotes the convex set over positive semi-definite matrices in $R^{n\times n}$. Suppose $X^{*}$ is an optimum of the problem with $rank(X^{*})\leq r$. The problem can be modeled as a non-convex problem, by writing $X=UU^{T}$, where U is an $n\times r$ matrix. Define $g(U)=f(UU^{T})$ and consider the direct optimization of the transformed problem:  \\
	$$ minimize_{U\in R^{n\times r}} g(U)\ where\ r\leq n.$$   
	Suppose $U^{*}$ is an optimum of the new problem. The method that uses gradient descent to solve the new optimization problem is called Factored Gradient Descent(FGD). By dropping convexity, [1] uses regularity condition to guarantee a linear convergence speed, similar to the case that the function is strongly convex and smooth. In this work, we apply adaptive gradient descent to minimize function $g(U)=f(UU^{T})$.\\
	~\\
	Such recasting of matrix optimization problems is widely popular, especially as the size of the matrix increases. This problem has many applications in solving machine learning problems.  The applications include matrix completion, affine rank minimization, matrix sensing, covariance selection, sparse PCA, phase retrieval and many other problems. \\
	~\\
	FGD ([1]) uses simple gradient descent to solve the optimization problem. Our goal in this note is to design a method which both performs well in the tests and keeps the theoretical guarantee. We apply adaptive gradient descent to solve the problem.  To the best of our knowledge, this is the first work that solves the re-parametrized problem with adaptive step size. Our main contributions can be summarized as follows: \\
	~\\
	1. We propose a simple algorithm to solve low-rank matrix factorization problem. The algorithm is computational efficient and it performs very well in all of the tests.\\
	~\\
	2. We analyze the local convergence for the algorithm and find that under the common condition for f and a good initialization, the convergence speed is better than FGD in all cases. For some conditions, the result can be much better than FGD.\\
	\subsection{Comparison with exisiting results}
	In literature, there are four main results for solving this problem. FGD[1], AGD-AC[2], AFGD[3] and Scaled-GD[4]. All these four methods use a fixed step size. Both AGD-AC and AFGD use Nesterov's accelerated gradient descent, but their initialization condition is very hard to satisfy (the proof of locally strongly convex requires the initialization point to be very close to the global minimum). For a wide range of problems that the condition number of f is not very large (for example, matrix factorization), the theoretical guarantee for these two methods is not better than FGD. Scaled-GD is a second order method that uses scaling matrices. The convergence speed of Scaled-GD is independent of the low-rank matrix, much faster than FGD. The theoretical guarantee of our method is not better than that of Scaled GD, but our method is worth studying because of its excellent performance in practice and its further use in more complex optimization problems, such as deep learning.
	\subsection{Roadmap}
	 The rest of this note is organized as follows. Section 2 contains basic notations. Section 3 presents the Adaptive Gradient Descent algorithm and the step size $\eta$ used and the reason for its selection. Section 4 contains the convergence guarantees of the algorithm. Section 5 contains some analysis of the theorem. Section 6 cointains the main proofs and lemmas of the main theorem and discusses some cases. Section 7 contains some tests. This note concludes with discussion on its furthur use in section 8.
	\section{Preliminaries}
	DIST and $R_{U}^{*}$ are defined the same as the [1], for matrices $U,V \in R^{n\times r}$,
	$$DIST(U,V)=min_{R:R\in O} \left\|U-VR\right\|_{F}.$$ O is the set of $r\times r$ orthonormal matrices R, such that $R^{T}R=I_{r\times r}$.  $R_{U}^{*}$ is the optimum R for U and $U^{*}$, which is $R_{U}^{*}=argmin_{R:R\in O}\|U-U^{*}R\|_{F}$. Let DIST refers to the true value and we use $\hat{DIST}$ as the estimation of DIST. For each epoch k, $\hat{DIST}(U_{k},U^{*})^{2}$ is $\hat{DIST}(U_{k},U^{*})^{2}=DIST(U_{k},U^{*})^{2}+\triangle_{k}$.  \\
	~\\
	We use $\|X\|_{F}$ and $\|X\|_{2}$ for the Frobenius and spectral norms of a matrix, respectively. Given a matrix X, we use $\sigma_{min}(X)$ and $\sigma_{max}(X)$ to denote the smallest and largest strictly positive singular values of X. If rank(X)=r, $\sigma_{r}(X)=\sigma_{min}(X)$. We  also use $\sigma_{1}(X)=\sigma_{max}(X)=\|X\|_{2}$. $Q_{A}$ denotes the basis of the column space of matrix A. For the function f, we only analyze the case that f is m-strongly convex and M-smooth. Define $\kappa =\frac{M}{m}$ as the condition number of f.
	\section{Adaptive Gradient Descent}
	I solve the non-convex problem via Adaptive Gradient Descent with update rule: 
	$$U_{k+1}=U_{k}-\eta_{k}\nabla f(U_{k}U_{k}^{T})\cdot U_{k}.$$
	$U_{k}$ is the value of U at epoch k and $\eta_{k}$ is the step size at epoch k. Adaptive gradient descent does this with a careful initialization (which is given below) and a special step size $\eta_{k}$. \\
	~\\
	This method is based on [1], which let $\eta_{k}=\eta=\frac{1}{16(M\left\|X_{0}\right\|_2+\left\| \nabla f(X_{0})\right\|_2)}$.  The convergence speed is $$DIST(U_{k+1},U^{*})^{2}\leq (1-\frac{3m\eta \sigma_{r}(X^{*})}{10})DIST(U_{k},U^{*})^{2}.$$ This can be proved by replacing $\eta_{k}$ with $\eta$ in (4) of section 6.\\
	~\\
	In my work, I use $\eta_{k}= \frac{4\hat{\eta}}{5}+\frac{3m\cdot \sigma_{r}(X^{*})\cdot \hat{DIST}(U_{k},U^{*})^{2}}{20\left\|\nabla f(X_{k})U_{k}\right\|_{F}^{2}}$, $\hat{\eta}$ will be defined in section 5, and $\hat{DIST}(U_{k},U^{*})^{2}$ is the estimation of $DIST(U_{k},U^{*})^{2}$ (In practice, we can replace $\frac{4\hat{\eta}}{5}$ with $\eta$ and replace  $\sigma_{r}(X^{*})$ with $\sigma_{r}(X_{0})$). According to the idea of [1], it uses the same step size for all the epochs, whcih ensures the descent inequality hold true for all the values of $\left\|\nabla f(X_{k})U_{k}\right\|_{F}^{2}$. This is not optimal, since we can increase the step size if the value of $\left\|\nabla f(X_{k})U_{k}\right\|_{F}^{2}$ is small. This leads us to think about using adaptive gradient descent to improve the upper bound of the descent inequality. Similar to [1], we use regularity condition to prove our theorem. The upper bound of $DIST(U_{k+1},U^{*})^{2}$ is a quadratic function of $\eta_{k}$ (see the right side of (4)). So we choose $\eta_{k}$ as the minimum of the quadratic function, which will ensure the least upper bound. This is the motivation of my work. \\
	~\\
	Then we will give the initialization condition.\\
	Initialization: Similar to [1], we assume that the method is initialized with a good starting point $X^{0}=U^{0}(U^{0})^{T}$ that is close to $X^{*}=U^{*}(U^{*})^{T}$. In particular, we assume
\begin{equation}
	DIST(U^{0},U^{*})\leq\frac{\sigma_{r}(U^{*})}{100\kappa}\frac{\sigma_{r}({X^{*}})}{\sigma_{1}({X^{*}})}.
\end{equation}
    This initialization condition may not be true in practice, nor is it necessary. But for theoretical proof, we assume it is true.
	\section{Theorem}
	In this part, we will give a theorem about the convergence speed of our work (adaptive grdient descent).
\begin{thm}
	Suppose $|\triangle_{k}| \leq \frac{DIST(U_{k},U^{*})^{2}}{2}$.  Then the inequality holds:
	$$DIST(U_{k+1},U^{*})^{2}\leq (1-\frac{3m\eta \sigma_{r}(X^{*})}{10})DIST(U_{k},U^{*})^{2}. $$ Also, the inequality holds,
	$$DIST(U_{k+1},U^{*})^{2}\leq (1-\frac{9m\eta^{*}_{k} \sigma_{r}(X^{*})}{80})DIST(U_{k},U^{*})^{2}.$$
\end{thm}
	\section{Analysis}
	In this part, we will analyze the result of the theorem. As we can know from the theorem, suppose the estimation is not very inaccurate, the result will not be worse than [1]. From the second inequality, we can know that the larger $\eta_{k}^{*}$ is, the faster the convergence speed will be. If $\frac{DIST(U_{k},U^{*})^{2}}{\left\|\nabla f(X_{k})U_{k}\right\|_{F}^{2}}$ is much larger than $\frac{\hat{\eta}}{m\sigma_{r}(X^{*})}$, $\eta_{k}^{*}$ is much larger than $\eta$, then the convergence speed is much faster than [1] at epoch k. Because the value of $\frac{\left\|\nabla f(X_{k})U_{k}\right\|_{F}^{2}}{DIST(U_{k},U^{*})^{2}}$ is between $\Theta (m^{2}\cdot \sigma_{r}(U^{*})^{4})$ and $\Theta (M^{2}\cdot \left\|U^{*}\right\|_{2}^{4})$, the bounds are too loose to be useful and the effect of adaptive gradient descent is not quite sure. If the situation that  $\eta_{k}^{*}$ is much larger than $\eta$ does not happen frequently, the improvement of convergence speed is limited. If the situation that  $\eta_{k}^{*}$ is much larger than $\eta$  happens frequently, there will be great improvement. 
	\section{Proof}
	In this section, I present the proof for the main theorem. First, there are some lemmas. \\
	Lemma1 (Lemma A.5 of the paper). For $X=UU^{T}$ that satisfies the assumption (1), let $\hat{\eta}=\frac{1}{16(M\left\|X\right\|_2+\left\| \nabla f(X)Q_{U}Q_{U}^{T}\right\|_2)}$, then $\hat{\eta}\geq \frac{5\eta}{6}$. \\
	~\\
	Lemma2 (Lemma 6.1 of the paper). For f being a M-smooth and (m,r)-strongly convex function, and under the assumption of (1), the  following inequality holds: 
\begin{equation}
	 <\nabla f(X)U, U-U^{*}R^{*}_{U}> \geq \frac{4\hat{\eta}}{5} \cdot \left\|\nabla f(X)U\right\|_{F}^{2}+\frac{3m}{20} \cdot \sigma_{r}(X^{*}) \cdot DIST(U,U^{*})^{2}.
\end{equation}
    ~\\
	Then because of $U_{k+1}=U_{k}-\eta_{k}\cdot \nabla f(X_{k})U_{k}$, the inequality holds.
\begin{equation}
\begin{split}
	DIST(U_{k+1},U^{*})^{2}=min_{R:R\in O} \left\|U_{k+1}-U^{*}R\right\|_{F}^{2} \leq \left\|U_{k+1}-U^{*}R_{U_{k}}^{*}\right\|_{F}^{2}\\
	= \left\|U_{k+1}-U_{k}\right\|_F^{2}+ \left\|U_{k}-U^{*}R_{U_{k}}^{*}\right\|_F^{2}-2<U_{k+1}-U_{k}, U^{*}R^{*}_{U_{k}}-U_{k}> \\ =\eta_{k}^{2}\cdot \left\|\nabla f(X_{k})U_{k}\right\|_F^{2}+ \left\|U_{k}-U^{*}R_{U_{k}}^{*}\right\|_F^{2} -2<U_{k+1}-U_{k},U^{*}R^{*}_{U_{k}}-U_{k}>.
\end{split}
\end{equation}
	Using lemma2, we can get.
\begin{equation}
\begin{split}
    DIST(U_{k+1},U^{*})^{2} \leq \eta_{k}^{2}\cdot \left\|\nabla f(X_{k})U_{k}\right\|_F^{2}+ \left\|U_{k}-U^{*}R_{U_{k}}^{*}\right\|_F^{2} \\ -2\eta_{k}(\frac{4\hat{\eta}}{5} \cdot \left\|\nabla f(X_{k})U_{k}\right\|_{F}^{2}+\frac{3m}{20} \cdot \sigma_{r}(X^{*}) \cdot DIST(U_{k},U^{*})^{2}).
\end{split}
\end{equation}
	It's obvious that $\left\|U_{k}-U^{*}R_{U_{k}}^{*}\right\|_F^{2}=DIST(U_{k},U^{*})^{2}$.
	Let $\eta_{k}^{*}$ be the value which minimizes the right side of the inequality, which is,
\begin{equation}
	\eta_{k}^{*}=\frac{4\hat{\eta}}{5}+\frac{3m\cdot \sigma_{r}(X^{*})\cdot DIST(U_{k},U^{*})^{2}}{20\left\|\nabla f(X_{k})U_{k}\right\|_{F}^{2}}.
\end{equation}
    ~\\
    Then we discuss these four cases, in section 5.1 and section 5.2, we assume we can get the precise value of $DIST(U_{k},U^{*})^{2}$, so we can use $\eta_{k}^{*}$ as the step size. In section 5.3 and section 5.4, we assume we can only get an estimation of $DIST(U_{k},U^{*})^{2}$, so we can only use the step size $\eta_{k}$, which lead to the conclusion of the theorem.
	\subsection{With step size $\eta_{k}^{*}$}
	Then with step size $\eta_{k}^{*}$,
\begin{equation}
	DIST(U_{k+1},U^{*})^{2} \leq DIST(U_{k},U^{*})^{2}
	-\frac{(\frac{3m}{10} \cdot \sigma_{r}(X^{*}) \cdot DIST(U_{k},U^{*})^{2}+\frac{8\hat{\eta}\cdot \left\|\nabla f(X_{k})U_{k}\right\|_{F}^{2}}{5})^{2}}{4\left\|\nabla f(X_{k})U_{k}\right\|_{F}^{2}}.
\end{equation}
	Using mean inequality, we have,
\begin{equation}
	DIST(U_{k+1},U^{*})^{2}\leq DIST(U_{k},U^{*})^{2}\cdot (1-\frac{12m\sigma_{r}(X^{*})\cdot \hat{\eta}}{25}).
\end{equation}
	Because of lemma1, $\hat{\eta}\geq \frac{5\eta}{6}$, the result is at least a bit better than the paper. 
	\subsection{$\eta_{k}^{*}$ is much larger than $\eta$}
	If $\frac{DIST(U_{k},U^{*})^{2}}{\left\|\nabla f(X_{k})U_{k}\right\|_{F}^{2}}$ is much larger than $\frac{\hat{\eta}}{m\sigma_{r}(X^{*})}$, this method can converge much faster than the paper, because of inequality(6) , 
\begin{equation}
	DIST(U_{k+1},U^{*})^{2}\leq DIST(U_{k},U^{*})^{2}(1-\frac{9m^{2}\sigma_{r}(X^{*})^{2}\cdot DIST(U_{k},U^{*})^{2}}{400\left\|\nabla f(X_{k})U_{k}\right\|_{F}^{2}}-\frac{6m\sigma_{r}(X^{*})\cdot \hat{\eta}}{25}).
\end{equation}
	According to the definition of $\eta_{k}^{*}$, we know that $\left\|\nabla f(X_{k})U_{k}\right\|_{F}^{2}=\frac{3m\cdot \sigma_{r}(X^{*})\cdot DIST(U_{k},U^{*})^{2}}{20(\eta_{k}^{*}-\frac{4\hat{\eta}}{{5}})}.$  \\
	Replacing $\left\|\nabla f(X_{k})U_{k}\right\|_{F}^{2}$ in (8), 
\begin{equation}
	DIST(U_{k+1},U^{*})^{2}\leq DIST(U_{k},U^{*})^{2}(1-\frac{3m\sigma_{r}(X^{*})\eta_{k}^{*}}{20}).
\end{equation}
	~\\
	The following two cases assume that we can only use the estimation of $DIST(U_{k},U^{*})^{2}$.
	\subsection{With step size $\eta_{k}$}
	Because of (4), we know that,
\begin{equation}
\begin{split}
    DIST(U_{k+1},U^{*})^{2} \leq DIST(U_{k},U^{*})^{2}+\left\|\nabla f(X_{k})U_{k}\right\|_F^{2} \cdot (\eta_{k}-\eta_{k}^{*})^{2} \\ -\left\|\nabla f(X_{k})U_{k}\right\|_F^{2}\cdot (\eta_{k}^{*})^{2}.
\end{split}
\end{equation}
    From the assumption of the theorem, we know that $|\eta_{k}-\eta_{k}^{*}|\leq \frac{\eta_{k}^{*}}{2}$. So,
\begin{equation}
    DIST(U_{k+1},U^{*})^{2} \leq DIST(U_{k},U^{*})^{2}-\frac{3\left\|\nabla f(X_{k})U_{k}\right\|_F^{2} \cdot (\eta_{k}^{*})^{2}}{4}.
\end{equation}
    Replacing $\eta_{k}^{*}$ with the precise value(5). We know that,
\begin{equation}
    DIST(U_{k+1},U^{*})^{2} \leq DIST(U_{k},U^{*})^{2}
    -\frac{3(\frac{3m}{10} \cdot \sigma_{r}(X^{*}) \cdot DIST(U_{k},U^{*})^{2}+\frac{8\hat{\eta}\cdot \left\|\nabla f(X_{k})U_{k}\right\|_{F}^{2}}{5})^{2}}{16\left\|\nabla f(X_{k})U_{k}\right\|_{F}^{2}}.
\end{equation}
    Using mean inequality, we have,
\begin{equation}
    DIST(U_{k+1},U^{*})^{2}\leq DIST(U_{k},U^{*})^{2}\cdot (1-\frac{9m\sigma_{r}(X^{*})\cdot \hat{\eta}}{25}).
\end{equation}
    Because of lemma1, $\hat{\eta}\geq \frac{5\eta}{6}$, the result is at least the same as [1]. 
    \subsection{$\eta_{k}^{*}$ is much larger than $\eta$}
	The most interesting case is that $\eta_{k}^{*}$ is much bigger than $\eta$, which leads to $\eta_{k}^{*}$ much larger than $\eta$. If $|\triangle_{k}|\leq \frac{DIST(U_{k},U^{*})^{2}}{2}$, then $\eta_{k}\geq \frac{\eta_{k}^{*}}{2}$. According to (12), 
\begin{equation}
DIST(U_{k+1},U^{*})^{2}\leq DIST(U_{k},U^{*})^{2}(1-\frac{27m^{2}\sigma_{r}(X^{*})^{2}\cdot DIST(U_{k},U^{*})^{2}}{1600\left\|\nabla f(X_{k})U_{k}\right\|_{F}^{2}}-\frac{9m\sigma_{r}(X^{*})\cdot \hat{\eta}}{50}).
\end{equation}
	Replacing $\left\|\nabla f(X_{k})U_{k}\right\|_{F}^{2}$, we can get,
\begin{equation}
    DIST(U_{k+1},U^{*})^{2}\leq (1-\frac{9m\eta^{*}_{k} \sigma_{r}(X^{*})}{80})DIST(U_{k},U^{*})^{2},
\end{equation}
    which is the second inequality of the theorem.
    \section{Tests}
    In the following matrix factorization tests, we compare adaptive gradient descent to the method of FGD. For convenience, we just use $\eta_{k}=\eta+\frac{3m\cdot \sigma_{r}(X_{0})\cdot DIST(U_{k},U^{*})^{2}}{20\left\|\nabla f(X_{k})U_{k}\right\|_{F}^{2}}.$ \\
    ~\\
    Matrix factorization problem is about:
    $$minimize_{ U\in R^{n\times r}} \left\|UU^{T}-A\right\|_{F}^{2}$$
    We design four tests to compare our method to FGD. In each test, each element of the matrix $U^{*}$ is taken from U(-1,1), and $U_{0}$ is taken either randomly near $U^{*}$ or randomly far from $U^{*}$. The relative error refers to $\frac{g(U_{k})}{g(U_{0})}$. The results of the tests are at the end of this note.  Figure 1 is a test that chooses U to be a $1000\times 2$ matrix. Figure 2 is a test that chooses U to be a $1000\times 5$ matrix. In each figure, there are two tests, one starting near the global minimum and the other starting far from the global minimum. From the tests, we can know that our method performs much better than FGD in all cases. Even if staring far from the global minimum, the additional term can help us achieve a faster convergence speed.
    \section{Conclusion}
    So, with the condition that $|\triangle_{k}|\leq \frac{DIST(U_{k},U^{*})^{2}}{2}$ and $\eta_{k}= \frac{4\hat{\eta}}{5}+\frac{3m\cdot \sigma_{r}(X^{*})\cdot \hat{DIST}(U_{k},U^{*})^{2}}{20\left\|\nabla f(X_{k})U_{k}\right\|_{F}^{2}}$, we can get a satisfying convergence speed. The convergence speed of [1] still holds for adaptive gradient descent. If the condition that $\eta_{k}^{*}$ is much bigger than $\eta$ happens frequently, the convergence speed will be much faster. However, the upper bound for $\left\|\nabla f(X_{k})U_{k}\right\|_{F}^{2}$ is too big to be useful. From the tests, we can know that the step size of FGD is very conservative, we can improve the empirical performance of gradient descent via using adaptive gradient descent and keep the guarantee of convergence speed at the same time.  The method can be applied to other optimization problems that use regularity condition to prove convergence.

    \begin{figure*}[h]
    	\centering
    	\includegraphics[scale=0.4]{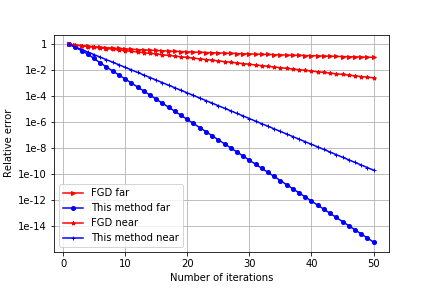}
    	\caption{n = 1000, r = 2}     \label{n = 1000, r = 2}
    \end{figure*}
    \begin{figure*}[h]
    	\centering
    	\includegraphics[scale=0.4]{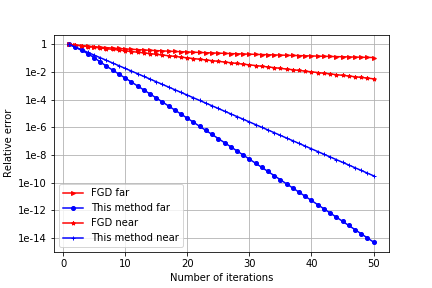}
    	\caption{n = 1000, r = 5}     \label{n = 1000, r = 5}
    \end{figure*}
\end{document}